\documentclass[11pt]{article}
\usepackage{amsfonts,amssymb,amsmath,graphicx}
\usepackage{color} % load the color package
   \definecolor{cites}{rgb}{0.50 , 0.00 , 0.00}  % colour for citations
   \definecolor{urls} {rgb}{0.00 , 0.00 , 0.50}  % colour for URL's
   \definecolor{links}{rgb}{0.00 , 0.00 , 0.50}   % colour for links
\usepackage[
      colorlinks=true,   % Yes, we want coloured links!
      citecolor=cites,   % Set the citations colour,
      urlcolor=urls,     % and the URL's colour,
      linkcolor=links,   % and the colour for the links.
      pdftitle={Two stable modifications of the finite section method},  % Put your title here
      pdfauthor={Marko Lindner}, % And your name here
      pdfpagemode=UseOutlines, % None, UseThumbs, UseOutlines, FullScreen
      pdfstartview=FitH,       % Fit, FitH, FitB
      bookmarksopen=false      % Bookmark tree already opened?
   ]{hyperref}

\newcommand\C{{\mathbb C}}
\newcommand\R{{\mathbb R}}
\newcommand\Z{{\mathbb Z}}
\newcommand\N{{\mathbb N}}
\newcommand\E{{\mathcal E}}
\newcommand\G{{\mathcal G}}
\newcommand\I{{\mathcal I}}

\renewcommand\H{{\mathcal H}}
\renewcommand\P{{\mathcal P}}

\newcommand\eps\varepsilon
\newcommand\BDO{{\rm BDO}}
\newcommand\diag{{\rm diag}}

\newcommand\conv{{\rm conv}}

\newcommand\im{{\rm im\,}}
\newcommand\ind{{\rm ind}}
\newcommand\sstab{\sigma^{\sf stab}}
\newcommand\opsp{\sigma^{\sf op}}

\newtheorem{theorem}{Theorem}[section]
\newtheorem{lemma}[theorem]{Lemma}

\newtheorem{proposition}[theorem]{Proposition}
\newtheorem{definition}[theorem]{Definition}

\newenvironment{remark}
 {\par\noindent\refstepcounter{theorem}{\bf Remark \thetheorem}\ }
 {\raisebox{1mm}{\framebox{}}\pagebreak[2]}

\newenvironment{example}
 {\par\noindent\refstepcounter{theorem}{\bf Example \thetheorem}\ }
 {\raisebox{1mm}{\framebox{}}\pagebreak[2]}

\newenvironment{proof}
 {\par\noindent{\bf Proof.}}
 {\rule{2mm}{2mm}\pagebreak[2]}

\begin{document}
\title{\bf Two stable modifications\\of the finite section method}
\author{{\sc Marko Lindner}}%\footnote{email: {\tt marko.lindner@mathematik.tu-chemnitz.de}}}
\date{\today}
\maketitle
\begin{quote}
\renewcommand{\baselinestretch}{1.0}
\footnotesize {\sc Abstract.} In this article we demonstrate and
compare two modified versions of the classical finite section method
for band-dominated operators in case the latter is not stable. For
both methods we give explicit criteria for their applicability.
\end{quote}


%%%%%%%%%%%%%%%%%%%%%%%%%%%%%%%%%%%%%%%%%%%%%%%%%%%%%%%%%%%%%%%%%%%%%%
\section{Introduction}
{\bf Infinite Matrices.} In this paper, we look at truncation
methods for the approximate solution of certain operator equations
$Au=b$ on the space $E:=\ell^p(\Z^N,X)$ of functions $u:\Z^N\to X$
with
\[
\|u\|\ =\ \left\{\begin{array}{cl}
\sqrt[\stackrel{p}{~}]{\sum\limits_{k\in\Z^N}
|u(k)|^p},&p\in[1,\infty)\\
\sup\limits_{k\in\Z^N}|u(k)|,&p=\infty
\end{array}\right\}\ <\ \infty,
\]
where $N\in\N$, $p\in[1,\infty]$ and $X$ is an arbitrary complex
Banach space. The operators $A$ that we have in mind are bounded
linear operators $E\to E$ which are induced, via
\begin{equation} \label{eq:induced}
(Au)(i)\ =\ \sum_{j\in\Z^N}a_{ij}\,u(j),\qquad i\in\Z^N,
\end{equation}
by a matrix $(a_{ij})_{i,j\in\Z^N}$ with operator entries
$a_{ij}:X\to X$. Among those operators we call $A$ a {\sl band
operator} if it is induced by a banded matrix, i.e. $a_{ij}=0$ if
$|i-j|$ is large enough, and we call $A$ a {\sl band-dominated
operator} and write $A\in\BDO(E)$ if $A$ is the limit, with respect
to the operator norm induced by the norm on $E$, of a sequence of
band operators. Also for $A\in\BDO(E)$, there is a unique (see
\cite[\S2.1.2]{RaRoSi:Book} or \cite[\S 1.3.5]{Li:Book}) matrix
$(a_{ij})_{i,j\in\Z^N}$ which induces $A$ via (\ref{eq:induced}); we
denote it by $[A]$.
\bigskip

{\bf Finite Sections.} If $A\in\BDO(E)$ is invertible then $Au=b$
has a unique solution $u\in E$ for every right-hand side $b\in E$.
An exact computation of $u$, however, is in general not possible
which is why one uses approximation methods. One of the most popular
approximation methods is as follows: Choose a sequence
$\Omega_1\subset\Omega_2\subset\cdots$ of finite subsets of $\Z^N$
that eventually covers every point of $\Z^N$ and replace the
infinite system
\[
Au\ =\ b\qquad\textrm{i.e.}\qquad \sum_{j\in \Z^N}a_{ij}\, u(j)\ =\
b(i),\quad i\in\Z^N
\]
by the sequence of finite systems
\begin{equation} \label{eq:FSM}
\sum_{j\in \Omega_n}a_{ij}\,\tilde u_n(j)\ =\ b(i),\quad
i\in\Omega_n
\end{equation}
for $n=1,2,...$\,. This procedure is called the {\sl finite section
method (FSM)}. The FSM is called {\sl applicable} if there exists an
$n_0\in\N$ such that, for every $b\in E$, (\ref{eq:FSM}) is uniquely
solvable for all $n\ge n_0$ and if the sequence $(\tilde u_n)$ of
solutions is bounded in $E$ and converges componentwise to the exact
solution $u$ of $Au=b$ as $n\to\infty$.

Here is how we will choose the finite sets $\Omega_1,\Omega_2,...$
in (\ref{eq:FSM}):

\begin{definition} \label{def:starlike}
We will say that $\Omega\subset\R^N$ is a {\sl valid starlike set}
if $\Omega$ is bounded, nonempty and has the property that, for
every $x\in\Omega$ and $\alpha\in [0,1)$, $\alpha x$ is an interior
point of $\Omega$.
\end{definition}

So in particular, $0$ is an interior point of every valid starlike
set. Moreover, all bounded convex sets $\Omega\subset\R^N$ with
interior point $0$ are valid starlike sets. Now, for every $n\in\N$,
put
\begin{equation} \label{eq:Omega}
\Omega_n\ :=\ n\Omega\,\cap\,\Z^N\qquad\textrm{and}\qquad P_n\ :=\
P_{\Omega_n},
\end{equation}
where, for a set $U\subseteq\Z^N$, by $P_U:E\to E$ we denote the
operator of multiplication by the characteristic function $\chi_U$
of $U$. Then we can abbreviate (\ref{eq:FSM}) as
\begin{equation} \label{eq:FSMshort}
P_n A P_n \tilde u_n\ =\ P_n b,\quad n=1,2,...\,.
\end{equation}
This truncation procedure is a very natural idea and the fact that
it can be performed on all infinite matrices creates the temptation
to simply use it and keep fingers crossed it will work. A positive
outcome, however, i.e. applicability as defined above, is in general
far from guaranteed. Here is the probably most elementary example
for which the FSM fails to apply:

\begin{example} \label{ex:shift}
Consider the shift operator $A=V_c:u\mapsto v$ on $E$ with
$u(k)=v(k+c)$ for every $k\in\Z^N$ and a fixed nonzero vector
$c\in\Z^N$. Then $V_c$ is invertible on $E$ but since $V_c$ maps
functions with support in $\Omega_n$ to functions supported in
$\Omega_n+c$, the truncated equation (\ref{eq:FSM}) alias
(\ref{eq:FSMshort}) is not solvable for general right-hand sides
(and even if it is solvable, the solution is not unique) -- no
matter how big $n$ is and how $\Omega$ is chosen.
\end{example}

Here is a slightly more sophisticated example:

\begin{example} \label{ex:blockdiag}
Let $N=1$ and consider the operator $A$ induced by the block
diagonal matrix
\[
\diag\left(\cdots,
\left(\begin{array}{cc}0&1\\1&0\end{array}\right),
\left(\begin{array}{cc}0&1\\1&0\end{array}\right), 1,
\left(\begin{array}{cc}0&1\\1&0\end{array}\right),
\left(\begin{array}{cc}0&1\\1&0\end{array}\right), \cdots\right)
\]
with the single $1$ entry at position zero. Then $A=A^{-1}$ is
invertible and, for $\Omega=[-1,1]$, its truncations $P_nAP_n$
correspond to the finite $(2n+1) \times (2n+1)$ matrices
\[
\diag\left(\left(\begin{array}{cc}0&1\\1&0\end{array}\right), \cdots
\left(\begin{array}{cc}0&1\\1&0\end{array}\right), 1,
\left(\begin{array}{cc}0&1\\1&0\end{array}\right), \cdots
\left(\begin{array}{cc}0&1\\1&0\end{array}\right)\right)
\]
if $n$ is even and to
\[
\diag\left(0, \left(\begin{array}{cc}0&1\\1&0\end{array}\right),
\cdots \left(\begin{array}{cc}0&1\\1&0\end{array}\right), 1,
\left(\begin{array}{cc}0&1\\1&0\end{array}\right), \cdots
\left(\begin{array}{cc}0&1\\1&0\end{array}\right), 0\right)
\]
if $n$ is odd. So the FSM (\ref{eq:FSMshort}) is not applicable
since all operators $P_nAP_n|_{\im P_n}$ with an odd $n$ are
non-invertible.
\end{example}

By \cite[Corollary 1.77]{Li:Book} (which is a consequence of
\cite[Theorem 6.1.3]{RaRoSi:Book}) one has that the FSM
(\ref{eq:FSMshort}) is applicable iff $A$ is invertible and the
sequence
\begin{equation} \label{eq:FSS}
(P_n AP_n+Q_n)_{n\in\N}
\end{equation}
is stable. Here we have put $Q_n:=I-P_n$ and we call a sequence
$(A_n)_{n\in\N}$ of operators $A_n:E\to E$ {\sl stable} if there
exists an $n_0\in\N$ such that all operators $A_n$ with $n\ge n_0$
are invertible and $\sup_{n\ge n_0} \|A_n^{-1}\|$ is finite. Also
note that $P_n AP_n+Q_n$ is invertible on $E$ iff $P_n AP_n$ is
invertible on the image of $P_n$ and that $\|(P_n
AP_n+Q_n)^{-1}\|=\max(1,\|(P_n AP_n|_{\im P_n})^{-1}\|)$.
\bigskip

The aim of this paper is to demonstrate two strategies, originally
developed in \cite{RaRoSi:FSMsubs,Li:FSMsubs} and
\cite{HeiLiPott:MSM}, that can be used if the FSM
(\ref{eq:FSMshort}) is not applicable:

{\bf Strategy 1: Pass to a subsequence. } As we have just seen, the
finite section method cannot be expected to work for every operator
$A$. But in some cases it is possible to ``adjust'' the method to
the operator at hand by choosing the right geometry $\Omega$ and an
appropriate subsequence of (\ref{eq:FSS}). The philosophy here is to
give the operator $A$ the chance to impose some of its
``personality'' on the (otherwise too ``impersonal'') method of
finite sections. In the previous example, for instance, one simply
has to remove all elements from the sequence (\ref{eq:FSS}) that
correspond to an odd value of $n$ to get a stable approximation
method for $A$ (or alternatively, one could replace $\Omega=[-1,1]$
by $[-2,2]$ and work with the whole sequence (\ref{eq:FSS})). We
believe that, for a given operator $A$, finding the right geometry
$\Omega$ and an appropriate sequence $n_1, n_2,...$ of natural
numbers such that the corresponding subsequence of finite sections
$P_{n_i}AP_{n_i}$ is stable (meaning that (\ref{eq:FSMshort}) is
only expected to be uniquely solvable, with solutions $\tilde u_n$
convergent to $u$, for a particular sequence $n=n_1,n_2,...$) is a
major task in the numerical analysis of the equation $Au=b$.

We will show that, under an additional condition on the operator
$A$, the finite section subsequence $(P_{n_i}AP_{n_i})_{i=1}^\infty$
is stable iff $A$ and every element from an associated set of
operators is invertible and the inverses are uniformly bounded. We
give a description of this associated set that depends on $A$,
$\Omega$ and the sequence $(n_i)_{i=1}^\infty$.

{\bf Strategy 2: Use rectangular instead of square systems. } As an
alternative approach to the FSM (\ref{eq:FSM}), we discuss the
slightly modified truncation scheme $P_m A P_n\,\tilde
u_{m,n}\approx P_m b$, i.e.
\begin{equation} \label{eqth:rFSM}
\sum_{j\in\Omega_n} a_{ij}\, \tilde u_{m,n}(j)\ \approx\ b(i),\qquad
i\in\Omega_m,
\end{equation}
leading to rectangular instead of quadratic finite subsystems of
$Au=b$ that are now to be solved approximately instead of exactly.

We prove that if $A$ is induced by a matrix $(a_{ij})$ with
$\|a_{ij}\|\to 0$ as $|i|\to\infty$ for every $j$ and if $A$ is
invertible then the modified method (\ref{eqth:rFSM}) is applicable.
By the latter we mean that, for every $\eps>0$ and every $b\in E$,
there exist $m_0,n_0\in\N$ and a precision $\delta>0$ such that all
(approximate) solutions of the rectangular system $\|P_m A P_n
\tilde u_{m,n}- P_m b\|<\delta$ with $m>m_0$ and $n>n_0$ are in the
$\eps$-neighbourhood of the exact solution $u$ of $Au=b$.

We also discuss how the two truncation parameters $m$ and $n$ are to
be coupled.
\bigskip

{\bf Short History.} The idea of the FSM is so natural that it is
difficult to give a historical starting point. First rigorous
treatments are from Baxter \cite{Baxter} and Gohberg \& Feldman
\cite{GohbergFeldman} on Wiener-Hopf and convolution operators in
dimension $N=1$ in the early 1960's. For convolution equations in
higher dimensions $N\ge 2$, the FSM goes back to Kozak \& Simonenko
\cite{Kozak,KozakSimonenko}, and for general band-dominated
operators with scalar \cite{RaRoSi1998} and operator-valued
\cite{RaRoSi2001,RaRoSi2001:FSM} coefficients, most results are due
to Rabinovich, Roch \& Silbermann. For the state of the art in the
scalar case for $p=2$, see \cite{Roch:FSM}.

The quest for stable subsequences if the FSM itself is instable is
getting more attention recently
\cite{RaRoSi:FSM_AP,RaRoSi:FSMsubs,SeidelSilbermann1,SeidelSilbermann2,Li:FSMsubs}.
In \cite{RaRoSi:FSMsubs}, the stability theorem for subsequences is
used to remove the uniform boundedness condition in dimension $N=1$.
Also the consideration of rectangular finite sections, although not
new in the numerical community, is now gaining more focus in the
numerical functional analysis literature (see
\cite{HeinigHellinger,Silb2003} for Toeplitz operators,
\cite{SeidelSilbermann1,SeidelSilbermann2} for band-dominated
operators and \cite{HeiLiPott:MSM} for even more general operators).
\bigskip

%%%%%%%%%%%%%%%%%%%%%%%%%%%%%%%%%%%%%%%%%%%%%%%%%%%%%%%%%%%%%%%%%%%%%%
\section{Strategy One: Stable Subsequences of the FSM}
\subsection{Preliminaries}
Let $E=\ell^p(\Z^N,X)$, $A\in\BDO(E)$ and $\Omega\subset\R^N$ be a
valid starlike set as in Definition \ref{def:starlike}. For an
infinite index set $\I=\{n_1,n_2,...\}\subseteq\N$, we study the
stability of the operator sequence
\begin{equation} \label{eq:FSsubs}
(P_nAP_n+Q_n)_{n\in\I}\ =\ (P_{n_i}AP_{n_i}+Q_{n_i})_{i=1}^\infty,
\end{equation}
where we suppose that $n_1, n_2,...$ is a strictly monotonous
enumeration of $\I$. For the study of this sequence as one item, we
will assemble it to a single operator. To do this, let
\begin{equation} \label{eq:Ai}
A_i\ :=\ \left\{\begin{array}{cl} P_{n_i}AP_{n_i}+Q_{n_i},&i\in\N,\\
I,&i\in\Z\setminus\N,\end{array}\right.
\end{equation}
put $E':=\ell^p(\Z^{N+1},X)$, thought of as $\ell^p(\Z,E)$, and
write $\oplus A_i$ for the map $u\mapsto v$ on $E'$ with
\begin{equation} \label{eq:stack}
%\left(\Big(\opl\limits_{i=1}^\infty A_i\Big)\ u\right)(j,k)
v(j,i)\ =\ \big(A_i\, u(\cdot,i)\big)(j), \qquad j\in\Z^N,\ i\in\Z.
\end{equation}
In other words, we think of $u\in E'$ as decomposed into layers
$u(\cdot,i)\in E$, $i\in\Z$, and let each $A_i$ act on the $i-$th
layer of $u$. We will therefore refer to $A_i$ as the $i-$th layer
of $\oplus A_i$. One can show that then $\oplus A_i\in\BDO(E')$.

A key argument in \cite{RaRoSi2001:FSM}, refined later in
\cite{Li:FSM,Li:Book,RaRoSi:Book,Roch:FSM}, is that the stability of
(\ref{eq:FSsubs}) is equivalent to $\oplus A_i$ being invertible at
infinity. Here we say that an operator $B\in\BDO(E')$ is {\sl
invertible at infinity} if there exist $C,D\in\BDO(E')$ and an
$m\in\N$ such that $CB\Theta_m=\Theta_m=\Theta_mBD$ holds, where
$\Theta_m$ is the operator of multiplication by the characteristic
function of $\Z^{N+1}\setminus\{-m,...,m\}^{N+1}$.

So it remains to study invertibility at infinity of $\oplus A_i$.
This is done in terms of so-called limit operators
\cite{RaRoSi:Book,CWLi:Memoir,Li:Book}. The idea is to reflect the
behaviour of an operator $B\in\BDO(E')$ at infinity by a family of
operators on $E'$ and to evaluate this family. To do this, we need
two notations. Firstly, for $B,B_1,B_2,...\in\BDO(E')$, we write
$B=\P'$-$\lim B_n$ if $[B_n]$ converges entrywise (in the norm of
$L(X)$) to $[B]$ as $n\to\infty$ and if $\sup_n\|B_n\|<\infty$.
Secondly, for $\alpha\in\Z^{N+1}$, let $V_\alpha':E'\to E'$ denote
the shift operator with $(V_\alpha' u)(k)=u(k-\alpha)$ for all
$k\in\Z^{N+1}$ and $u\in E'$.

If $B\in\BDO(E')$, $h=(h(1),h(2),...)\subseteq\Z^{N+1}$ is a
sequence with $|h(n)|\to\infty$ and the operator sequence
$V_{-h(n)}'BV_{h(n)}'$ is $\P'-$convergent as $n\to\infty$ then its
limit will be denoted by $B_h$ and is called {\sl limit operator of
$B$ w.r.t. the sequence $h$}. In an analogous fashion, one defines
limit operators in $\BDO(E)$. To distinguish between operators on
$E'$ and on $E$ we write $\P$-$\lim$ and $V_\alpha$ with
$\alpha\in\Z^N$ if we are in the $E$ setting. Different sequences
$h$ generally lead to different limit operators and often the
sequence $V_{-h(n)}BV_{h(n)}$ does not $\P-$converge at all. We will
call $B\in\BDO(E)$ a {\sl rich} operator if every sequence
$h=(h(1),h(2),...)\subseteq\Z^N$ with $|h(n)|\to\infty$ has a
subsequence $g$ such that the limit operator $B_g$ exists.

As a final preparation, we turn our attention to the geometry of
$\Omega$. Let $\Gamma:=\partial\Omega$ be the boundary of $\Omega$
and, for every $n\in\N$, put
\[
\Gamma_n\ :=\ (n\Gamma\,+\,H)\,\cap\,\Z^N\qquad\textrm{with}\qquad
H=\left(-1/2\,,\, 1/2\,\right]^N
\]
and then let
\[
\Gamma_\I\ :=\ \bigcup_{n\in\I}\Gamma_n.
\]
For a sequence $h=(h(1),h(2),...)\subseteq \Gamma_\I$, say
$h(k)\in\Gamma_{m_k}$ for some $m_k\in\I$, and a set
$S\subseteq\Z^N$, we call $S$ the {\sl geometric limit of $\Omega$
w.r.t. $h$} and write $S=\Omega_h$ if, for every $m\in\N$, there
exists a $k_0\in\N$ such that
\[
\big(\Omega_{m_k}-h(k)\big)\,\cap\,\{-m,...,m\}^N\ =\
S\,\cap\,\{-m,...,m\}^N,\qquad k\ge k_0.
\]
Note that in this case $V_{-h(k)}P_{m_k}V_{h(k)}$ is $\P-$convergent
to $P_S$ as $k\to\infty$. For a polytope $\Omega$, the only
candidates for the geometric limit $S$ w.r.t a sequence $h\subseteq
\Gamma_\I$ are intersections of finitely many half spaces and $\Z^N$
(discrete half spaces, edges, corners, etc.).

\subsection{The Stability Theorem for Subsequences}
Given a rich operator $A\in\BDO(E)$ on $E=\ell^p(\Z^N,X)$ with
$p\in[1,\infty]$, $N\in\N$ and a complex Banach space $X$, a valid
starlike set $\Omega\in\R^N$, and an index set
$\I=\{n_1,n_2,...\}\subseteq\N$ with $n_1<n_2<\cdots$, we put
\[
\H_{\Omega,\I}(A)\ :=\ \big\{\,h=(h(1),h(2),...)\,:\,
h(k)\in\Gamma_\I\ \forall k,\, |h(k)|\to\infty,\, A_h \textrm{
exists, } \Omega_h \textrm{ exists}\,\big\}
\]
and
\begin{equation} \label{eq:sstab}
% \opsp_{\Omega,\I}(A)\ :=\ \big\{A_h\,:\,h\in\H_{\Omega,\I}(A)\big\},\quad
\sstab_{\Omega,\I}(A)\ :=\
\{A\}\cup\big\{P_{\Omega_h}A_hP_{\Omega_h}+Q_{\Omega_h}\,:\,h\in\H_{\Omega,\I}(A)\big\}.
\end{equation}
Then the following theorem holds.

\begin{theorem} \label{th:main}
Under the conditions mentioned above, the following are equivalent.

\begin{tabular}{rl}
(i)& The sequence $(P_{n_i}AP_{n_i}+Q_{n_i})_{i=1}^\infty$ is
stable.\\[1mm]
(ii)& The operator $\oplus A_i$, with $A_i$ as in (\ref{eq:Ai}), is invertible at infinity.\\[1mm]
~(iii)& All operators in $\sstab_{\Omega,\I}(A)$ are invertible with
their inverses uniformly bounded.
\end{tabular}
\end{theorem}

\begin{proof}
See Theorems 3.1 and 3.5 in \cite{Li:FSMsubs}.
\end{proof}

For dimension $N=1$, our statement coincides with a two-sided
version of \cite[Theorem 3]{RaRoSi:FSMsubs}. As such it generalizes
\cite[Theorem 3]{RaRoSi2001:FSM} (also see \cite[Theorem
6.2.2]{RaRoSi:Book}, \cite[Theorem 4.2]{Li:Book} and \cite[Theorem
2.7]{Roch:FSM}) from the full sequence $\I=\N$ to arbitrary infinite
subsequences with index set $\I\subseteq\N$. For $N=2$ and $\Omega$
a convex polygon with integer vertices, our Theorem \ref{th:main},
together with (\ref{eq:sstab}), corrects another version of the
stability spectrum %(see (16) and Example 4.1 in \cite{Li:FSM})
(see (\ref{eq:RaRoSisstab}) and Example \ref{ex:RaRoSiError} below)
that was previously suggested in the literature (see
\cite{RaRoSi2001:FSM,RaRoSi:Book}) for $\I=\N$. Moreover, our result
demonstrates how to deal with subsequences $\I\subseteq\N$ by
restricting consideration to sequences $h=(h(1),h(2),...)$ with
values in the set $\Gamma_\I=\cup_{n\in\I}\Gamma_n$. For dimensions
$N>2$, to our knowledge, the result is new -- even in cases like
$\I=\N$ or $\Omega$ a convex polytope.

\subsection{Examples}
As a particularly illustrative and not too difficult class of
examples, we will look at operators that are induced by an adjacency
matrix. Therefore, put $X=\C$, let $\E$ denote a set of pairwise
disjoint doubletons $\{i,j\}$ (i.e. sets  $\{i,j\}=\{j,i\}$ with
exactly two elements) with $i,j\in\Z^N$, $i\ne j$, and put
\[
a_{ij}\ :=\ \left\{\begin{array}{cl} 1,&\textrm{if }\{i,j\}\in\E\textrm{ or } i=j\not\in\bigcup\limits_{e\in\E} e,\\
0,&\textrm{otherwise,}\end{array}\right.
\]
for all $i,j\in\Z^N$. Then $(a_{ij})_{i,j\in\Z^N}$ is the {\sl
extended adjacency matrix} of the undirected graph $\G=(\Z^N,\E)$
with vertex set $\Z^N$ and edges $\E$. We write $Adj(\G)$ for the
operator that is induced by this matrix $(a_{ij})$ and note that
$Adj(\G)$ is band-dominated iff $b:=\sup_{\{i,j\}\in\E}|i-j|$ is
finite, in which case $Adj(\G)$ is even a band operator with
band-width $b$.

If applied to an element $u\in E=\ell^p(\Z^N,X)$, the operator
$Adj(\G)$ ``swaps'' the values $u(i)$ and $u(j)$ around if $\{i,j\}$
is an edge of $\G$, and it leaves all values $u(k)$ untouched for
which $k\in\Z^N$ is not part of an edge of $\G$. From this it is
obvious that $\|Adj(\G)\|=1$ and that $Adj(\G)$ is invertible and
coincides with its inverse. Moreover, it is clear that, for
$n\in\N$, the $n$-th finite section $P_nAdj(\G) P_n+Q_n$ is
invertible iff each edge $e\in\E$ has either both or no vertices in
$\Omega_n=n\Omega\cap\Z^N$. In the latter case, $P_nAdj(\G) P_n+Q_n$
equals $Adj(\G_n)$, where $\G_n=(\Z^N,\E\cap\Omega_n^2)$, is again
its own inverse and has norm 1. So we get that, for $A=Adj(\G)$, the
sequence (\ref{eq:FSsubs}) is stable iff, for all sufficiently large
$n\in\I$, each edge $e\in\E$ has either both or no vertices in
$\Omega_n$.

Note that Example \ref{ex:blockdiag} was already of the form
$A=Adj(\G)$, namely with $N=1$ and
\[
\E\ =\ \Big\{\,...,\{-4,-3\},\{-2,-1\},\{1,2\},\{3,4\},...\,\Big\}.
\]
Here $\Omega_n$ separates the vertices of the edge $\{-n-1,-n\}$ and
also of $\{n,n+1\}$ if $n$ is odd.

We continue with two examples demonstrating that two particular sets
of operators that are closely related to $\sstab_\Omega(A)$ -- and
that have, in the past, been suggested to replace (\ref{eq:sstab})
in the $N=2$, $\I=\N$ version of Theorem \ref{th:main} -- are
actually not stability spectra (meaning that Theorem \ref{th:main}
is incorrect for $\I=\N$ with $\sstab_\Omega(A)$ replaced by any of
them) if $N>1$. These two ``non-replacements'' for
$\sstab_\Omega(A)$ are
\begin{equation} \label{eq:RaRoSisstab}
\{A\}\ \cup\ \bigcup_{x\in\Gamma}
\{P_{\Omega_x}BP_{\Omega_x}+Q_{\Omega_x}\ :\ B\in\opsp_{x}(A)\}
\end{equation}
and
\begin{equation} \label{eq:Sisstab}
\{A\}\ \cup\ \bigcup_{x\in\Gamma}
\{P_{\Omega_x}BP_{\Omega_x}+Q_{\Omega_x}\ :\ B\in\opsp_{x,{\sf
ray}}(A)\},
\end{equation}
where $\Gamma=\partial\Omega$ and, for every $x\in\Gamma$,
$\Omega_x\subseteq\Z^N$ is the limit of $n(\Omega-x)\cap\Z^N$ as
$n\to\infty$ in the sense that, for each $m\in\N$,
\[
n(\Omega-x)\cap\{-m,...,m\}^N\ =\ \Omega_x\cap\{-m,...,m\}^N
\]
for all sufficiently large $n\in\N$. Finally, $\opsp_{x}(A)$ is the
set of all limit operators $A_h$ of $A$ with respect to sequences
$h=(h(1),h(2),...)\subseteq\Z^N$ going to infinity in the direction
$x$, i.e. $h(n)/|h(n)|\to x/|x|$, and $\opsp_{x,{\sf ray}}(A)$ is
the set of all limit operators $A_h$ with respect to sequences of
the form $h=([m_1\,x],[m_2\,x],...)\subseteq\Z^N$ where $(m_n)$ is
an unbounded monotonously increasing sequence of positive reals and
$[\,\cdot\,]$ means componentwise rounding to the nearest integer.

\begin{example} \label{ex:RaRoSiError}
Take $N=2$, $\Omega=[-1,1]^2$ and let $A=Adj(\G)$ with
$\G=(\Z^2,\E)$ and
\[
\E=\Big\{\{\,(k^2-k-1,k^2)\,,\,(k^2-k,k^2)\,\}\ :\ k=1,2,...\Big\}.
\]
Then, with respect to $h=(h(1),h(2),...)$ with
$h(k)=(k^2-k-1,k^2)\in\Z^2$, the limit operator of $A$ exists and is
equal to $B=Adj(\G')$, where
$\G'=\Big(\Z^2,\Big\{\{(0,0),(1,0)\}\Big\}\Big)$. Since
$h(k)/|h(k)|\to x/|x|$ with $x=(1,1)$, we have that
$B\in\opsp_x(A)$. But $\Omega_x=\{...,-1,0\}^2$ separates $(0,0)$
from $(1,0)$ so that
$P_{\Omega_x}BP_{\Omega_x}+Q_{\Omega_x}\in(\ref{eq:RaRoSisstab})$ is
not invertible. However, the whole finite section sequence
(\ref{eq:FSS}) is stable since all edges $e\in\E$ have either both
or no points in $\Omega_n$, so that $P_n AP_n+Q_n=Adj(\G_n)$ with
$\G_n=(\Z^2,\E\cap\Omega_n^2)$ for every $n\in\N$. So
(\ref{eq:RaRoSisstab}) is not a valid replacement of
(\ref{eq:sstab}) as stability spectrum.

Note that the element of (\ref{eq:sstab}) that corresponds to the
limit operator $B=A_h$ of $A$ is
$P_{\Omega_h}BP_{\Omega_h}+Q_{\Omega_h}$ with
$\Omega_h=\Z\times\{...,-1,0\}$ instead of $\{...,-1,0\}^2$, which
is again equal to $B$ (since both $(0,0)$ and $(1,0)$ are in
$\Omega_h$) and hence invertible.
\end{example}

Similarly, we can rule out (\ref{eq:Sisstab}) as stability spectrum
by the following example:

\begin{example} \label{ex:SiError}
Again take $N=2$, $\Omega=[-1,1]^2$ and let $A=Adj(\G)$ with
$\G=(\Z^2,\E)$ and
\[
\E=\Big\{\{\,(k^2-k,k^2)\,,\,(k^2-k,k^2+1)\,\}\ :\ k=1,2,...\Big\}.
\]
Then, with respect to $h=(h(1),h(2),...)$ with
$h(k)=(k^2-k,k^2)\in\Z^2$, the limit operator of $A$ exists and is
equal to $B=Adj(\G')$, where
$\G'=\Big(\Z^N,\Big\{\{(0,0),(0,1)\}\Big\}\Big)$.

Again $B\in\opsp_x(A)$ with $x=(1,1)$. But $B\not\in\opsp_{x,{\sf
ray}}(A)$ neither is $B$ in $\opsp_{y,{\sf ray}}(A)$ for any other
$y\in\Gamma$! In fact, it holds that $\opsp_{y,{\sf ray}}(A)=\{I\}$
for all $y\in\Gamma$, whence (\ref{eq:Sisstab}) is elementwise
invertible with uniformly bounded inverses. However, the finite
section sequence (\ref{eq:FSS}) is not stable since $\Omega_n$
separates $(k^2-k,k^2)$ from $(k^2-k,k^2+1)$ if $n=k^2$. So also
(\ref{eq:Sisstab}) is not a valid replacement of (\ref{eq:sstab}) as
stability spectrum.

Note that, for $\I=\N$, (\ref{eq:sstab}) contains the operator
$P_{\Omega_h}BP_{\Omega_h}+Q_{\Omega_h}$ with
$\Omega_h=\Z\times\{...,-1,0\}$, which is non-invertible since
$\Omega_h$ separates $(0,0)$ from $(0,1)$. This operator is however
removed from (\ref{eq:sstab}) if we remove all (sufficiently large)
square numbers from $\I$, which matches our observation that
$P_nAP_n+Q_n$ is non-invertible iff $n$ is a square.
\end{example}

It is clear that Examples \ref{ex:RaRoSiError} and \ref{ex:SiError}
can easily be heaved to dimensions $N>2$. Let us look at another
example, for simplicity also in dimension $N=2$.

\begin{example} \label{ex:diamond}
We look at $A=Adj(\G)$ for $\G=(\Z^2,\E)$, where
\[
\E\ =\ \Big\{ \{\,(k,1)\, ,\, (k+1,0)\,\}\ :\ k=1,2,... \Big\}.
\]
It is not hard to see that every limit operator of $A$ is either the
identity operator $I$ or the operator $B=Adj(\G')$ for
$\G'=(\Z^2,\E')$, where
\[
\E'\ =\ \Big\{ \{\,(k,1)\, ,\, (k+1,0)\,\}\ :\ k\in\Z \Big\},
\]
or it is a translate of $B$. Looking at $B$ and noting that $B=A_h$
for all sequences $h=(h(1),h(2),...)$ with $h(k)=(m_k,0)$ and
$m_k\to +\infty$, we can say how $\Omega$ has to look locally at the
intersection $z$ of its boundary $\Gamma$ with the positive $x$-axis
in order for the finite section method to be stable. Here the upward
tangent of $\Gamma$ at $z$ has to enclose an angle $\alpha\in
(90^o,135^o]$ with the positively directed $x$-axis. So, for
example, the finite section sequence is stable if $\Omega$ is the
square $\conv\{(1,0),(0,1),(-1,0),(0,-1)\}$ or the triangle
$\conv\{(0,2),(2,-2),(-2,-2)\}$, whereas it does not even have a
stable subsequence if $\Omega$ is the square $[-1,1]^2$.
\end{example}

The next example is closely related to Example \ref{ex:blockdiag}.

\begin{example} \label{ex:noway}
{\bf a) }Let $A=Adj(\G)$ where $\G=(\Z,\E)$ is the following infinite graph:\\
\includegraphics[width=\textwidth]{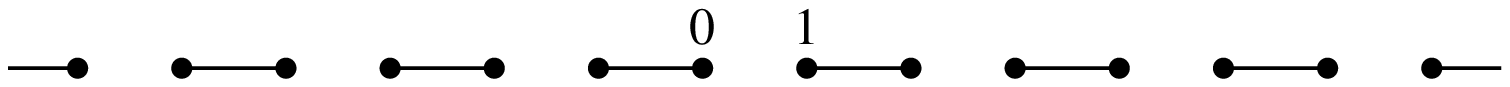}\\
Then, no matter how we choose $\Omega=[a,b]$ with integers $a<0<b$,
the finite section method does not even have a stable subsequence. A
workaround would be to take $\Omega=[-1,1)$ or to increase the
dimension to $N=2$, where we place the edges $\E$ along the $x$-axis
and put $\Omega=\conv\{(-1,0),(1,1),(0,-1)\}$, for example. In the
latter case, the finite section subsequence corresponding to
$\I=4\N+1$ turns out to be stable.

{\bf b) } In contrast to a), there is no workaround whatsoever if
$A=Adj(\G)$ with the following graph $\G$ (embedded in dimension
$N=1$ or higher):\\
\includegraphics[width=\textwidth]{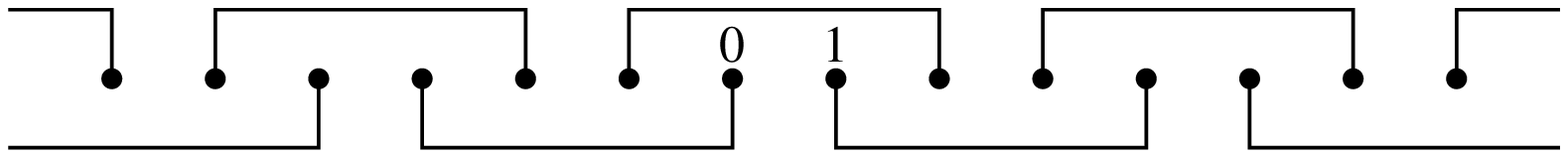}
For every valid starlike set $\Omega$ and every $n\in\N$, the set
$\Omega_n$ separates the endpoints of at least two edges of $\G$ so
that $P_nAP_n+Q_n$ is non-invertible.
\end{example}

For any dimension $N\in\N$, any valid set $\Omega\in\R^N$ and any
given sequence $n_1<n_2<\cdots$ of naturals, one can construct a
graph $\G$ in the style\footnote{The idea is to take the graph from
Example \ref{ex:noway} b) and to place ``gaps'' between $a_i:=\lceil
a n_i\rceil$ and $a_i-1$ and between $b_i:=\lfloor b n_i\rfloor$ and
$b_i+1$ for $i=1,2,...$, where $a<0$ and $b>0$ are the unique
intersection points of $\Gamma=\partial\Omega$ with the $x$-axis and
$\lceil \cdot\rceil$ and $\lfloor \cdot\rfloor$ stand for rounding
up and down to the next integer, respectively.} of Example
\ref{ex:noway} b)  such that $(P_nAP_n+Q_n)_{n\in\I}$ is stable iff
$\I$ is a subset of $\{n_1,n_2,...\}$.

\subsection{Some Words on the Case $N=1$}
Not surprisingly, the results are most complete in dimension $N=1$,
where one can sharpen and extend much of what was said previously
(also see
\cite{RaRoSi:FSMsubs,RaRoSi:IndexFSM,SeidelSilbermann1,SeidelSilbermann2,Li:FSMsubs}).
This is clearly due to the simple geometry of this setting: Firstly,
to infinity there are only two ways to go: right or left, and
secondly, all valid starlike sets are intervals from $a$ to $b$ with
reals $a<0<b$ so that there are only two possibilities for
$\Omega_h$ in (\ref{eq:sstab}): $\{0,1,...\}$ and $\{...,-1,0\}$.

The main result on $N=1$ is by Rabinovich, Roch and Silbermann. It
highlights an important benefit from extending the stability theorem
from the full finite section sequence to subsequences. A proof can
be found in \cite{RaRoSi:FSMsubs} or, slightly generalized, in
\cite[\S6]{Li:FSMsubs}.

\begin{proposition} \cite{RaRoSi:FSMsubs}
The uniform boundedness condition in Theorem \ref{th:main} (iii) is
redundant if $N=1$.
\end{proposition}

As a next result in dimension $N=1$, we mention that if the full FSM
(\ref{eq:FSS}) is stable for one valid $\Omega$ then (\ref{eq:FSS})
has a stable subsequence for all valid $\Omega$. So conversely, if
there exists a valid $\Omega$ for which (\ref{eq:FSS}) has no stable
subsequence then there is no valid $\Omega$ for which the whole
sequence (\ref{eq:FSS}) is stable. A proof can be found in
\cite[\S5]{Li:FSMsubs}.

Example \ref{ex:noway} b) has shown that, for some operators, the
finite section method cannot be ``adjusted'', via choosing $\Omega$
and $\I$, to become stable. We now give a necessary criterion for
the existence of an index set $\I\subseteq\N$ and a valid $\Omega$
such that (\ref{eq:FSsubs}) is stable.

\begin{proposition} \label{prop:N=1:I}
Let $E=\ell^p(\Z,\C)$ with $p\in[1,\infty]$ and $A\in\BDO(E)$. For
the existence of a valid starlike set $\Omega\subset\R$ and an
infinite index set $\I\subseteq\N$ such that the sequence
$(P_nAP_n+Q_n)_{n\in\I}$ is stable it is necessary, but not
sufficient, that $A$ is invertible and the Fredholm index
$\ind_+(A):=\ind(P_\N AP_\N|_{\im P_\N})$ is zero.
\end{proposition}

A proof, based on work on the Fredholm index of band-dominated
operators in \cite{RaRoRoe}, can be found in \cite[\S5]{Li:FSMsubs}.
In \cite{LiRaRo} (see \cite{Roch:ellp,Li:Wiener} for $p\ne 2$) we
have shown that, under the additional condition that all diagonals
of $[A]$ are slowly oscillating, invertibility of $A$ and
$\ind_+(A)=0$ are even sufficient for the stability of the full
finite section sequence (\ref{eq:FSS}) for all valid $\Omega$. Here
we call a sequence $(b_k)_{k\in\Z}$ {\sl slowly oscillating} if
$b_{k+1}-b_k\to 0$ as $k\to\pm\infty$.

\begin{remark} \label{rem:Gohberg}
By Proposition \ref{prop:N=1:I}, for an invertible operator $A$ with
$\kappa:=\ind_+(A)\ne 0$, there is no valid $\Omega$ and no index
set $\I\subseteq\N$ for which (\ref{eq:FSsubs}) is stable. This
problem of a nonzero plus-index $\kappa$ can be overcome as follows:
Instead of solving $Au=b$, one looks at $V_\kappa Au=V_\kappa b$
with $V_\kappa$ as in Example \ref{ex:shift}. Since $V_\kappa$ is
invertible, these two equations are equivalent. Moreover, also
$A':=V_\kappa A$ is invertible and
\[
\ind_+(A')\ =\ \ind_+(V_\kappa A)\ =\ \ind_+(V_\kappa)+\ind_+(A)\ =\
-\kappa + \kappa\ =\ 0.
\]
This preconditioning-type procedure of shifting the whole system
(all matrix entries and the right hand side $b$) down by $\kappa$
rows is reminiscent of Gohberg's statement that, in a two-sided
infinite matrix, ``it is every diagonal's right to claim to be the
main one'' (see page 51f in \cite{HaRoSi:2001}). Our results show
that, however, from the perspective of the FSM, there is one
diagonal that deserves being the main diagonal a bit more than the
others.
\end{remark}

%%%%%%%%%%%%%%%%%%%%%%%%%%%%%%%%%%%%%%%%%%%%%%%%%%%%%%%%%%%%%%%%%%%%%%
\section{Strategy Two: Rectangular Subsystems}
\subsection{Motivation}
Let us go back to Example \ref{ex:shift} and try to fix one of the
basic problems of the FSM. For simplicity, think of dimension $N=1$
and $\Omega=[-1,1]$ so that $\Omega_n=\{-n,...,n\}$ for all
$n\in\N$. Look at the shift operator $A=V_k$ with $k=1$, say. The
FSM for the solution of $Au=b$, that is
\begin{equation} \label{eq:FSM4MSM}
P_nAP_n u_n\ =\ P_n b,\qquad n=1,2,...,
\end{equation}
thinks of an approximate solution $u_n$ with support in
$\{-n,...,n\}$, then applies the operator -- in our case the forward
shift by $1$ component -- and afterwards cuts off at $\{-n,...,n\}$
again, hereby trying to match the restriction of the right-hand side
$b$ to $\{-n,...,n\}$. It is clear that this truncated equation
(\ref{eq:FSM4MSM}) is in general not solvable since the left-hand
side of (\ref{eq:FSM4MSM}) always has a $0$ at component $-n$
whereas the right-hand side has the same component $-n$ as $b$ has.
Even if $b(-n)=0$ and (\ref{eq:FSM4MSM}) is solvable then the
solution $u_n$ is not unique\footnote{Of course, for a finite
quadratic system, solvability for all right-hand sides (i.e.
surjectivity of the finite matrix operator) is equivalent to
uniqueness of the solution (i.e. injectivity). The approach here is
to say that lack of surjectivity can be overcome by looking for
approximate rather than exact solutions, whereas lack of injectivity
is a more serious problem that will be dealt with by adding more
equations (i.e. more matrix rows) to the finite system.} since its
$n$th component got shifted and then cut off whence it is irrelevant
for (\ref{eq:FSM4MSM}).

The observation generalizes to band and band-dominated operators of
course. If we truncate $u_n$ at $\{-n,...,n\}$ and apply a band
operator $A$ with band-width $w$ then $AP_nu_n$ is supported in
$\{-n-w,...,n+w\}$ whence, for the same reasons as illustrated for
the shift $A=V_1$, it is better to cut off at $\{-m,...,m\}$ with
$m=n+w$ and not at $\{-n,...,n\}$. The resulting system
\begin{equation} \label{eq:MSM'}
P_mAP_n u_n\ =\ P_m b,\qquad n=1,2,...,
\end{equation}
with $m=n+w$ is over-determined -- it has rectangular matrices that
have $2w$ more rows than they have columns. But one can still try to
solve it approximately (by least squares, say).

From the matrix point of view, $[AP_n]$ is the same as $[A]$, only
with all but columns number $-n,...,n$ put to zero. If the
horizontal cut-off $[P_mAP_n]$ (that one also has to do to get a
finite system for the computer) is done at $m=n$, like in
(\ref{eq:FSM4MSM}), then some `large' entries of $[AP_n]$ will get
cut off (recall $A=V_1$) which might cause problems as mentioned
earlier; so it could be good to choose $m$ a bit larger. In fact, if
$A$ has the property
\begin{equation} \label{eq:MSMdecay}
\|P_mAP_n\ - AP_n\|\ \to\ 0,\qquad \textrm{i.e.}\qquad \|Q_mAP_n\|\
\to\ 0\qquad\textrm{as}\qquad m\to\infty
\end{equation}
for all $n\in\N$ then it seems possible to work with this
rectangular cut-off idea, where $m$ in (\ref{eq:MSM'}), depending on
$n$, is chosen large enough to make $\|P_mAP_n-AP_n\|=\|Q_mAP_n\|$
small enough. The class of operators with property
(\ref{eq:MSMdecay}) clearly contains all of $\BDO(E)$.

The above idea is so natural that it can hardly be new. Indeed, it
is already used by some of the numerical community and it goes back
at least to the 1960's when Cleve Moler suggested, roughly speaking:
If square submatrices give you problems, make them higher and use
least squares. In \cite{HeiLiPott:MSM} we have not only reinvented
this method, we have (and that seems to be new) given a proof that,
in the setting of a rather general Banach space $E$, the method is
applicable as soon as $A$ is invertible and subject to
(\ref{eq:MSMdecay}). We will now recall the main steps of this
proof.

\subsection{The Rectangular Finite Section Method (rFSM)}
We will work with the same spaces $E=\ell^p(\Z^N,X)$ and the same
projection operators $P_n$ and $Q_n$ here as defined above, but we
will now exclude the case $p=\infty$ because we require strong
convergence $P_n\to I$, i.e. $P_nu\to u$ for all $u\in E$, as
$n\to\infty$. The subspace
$c_0(\Z^N,X)=\{u\in\ell^\infty(\Z^N,X):P_n u\to u\}$ of
$\ell^\infty(\Z^N,X)$ is however a valid choice for $E$.

Now suppose $A:E\to E$ is a bounded and invertible linear operator
with (\ref{eq:MSMdecay}), that means $\|a_{ij}\| \to 0$ as
$|i|\to\infty$ for every fixed $j\in\Z^N$, where $[A]=(a_{ij})$ is
the matrix representation of $A$. Then the equation $Au=b$ has a
unique solution $u=:u_0$ for every right-hand side $b\in E$. For the
approximate computation of $u_0$ we propose the following method:
For given precision $\delta>0$ and cut-off parameters $m$ and
$n\in\N$, calculate a solution $u\in \im P_n$ of the inequality
\begin{equation} \tag{rFSM}%\label{eq:rFSM}
\|P_{m} A P_{n} u - P_{m}b\| \ <\  \delta.
\end{equation}

We start with a result about the existence of solutions of (rFSM).

%-----------------------------------------------------
% Existence Theorem
%-----------------------------------------------------
\begin{definition}
We say that $n_0\in\N$ is an {\sl admissible $n$-bound} for $A$, $b$
and a given precision $\delta>0$ if {\rm (rFSM)} is solvable in $E$
for all $m\in\N$ and all $n\ge n_0$.
\end{definition}

\begin{proposition} \label{prop:exist}
For every $\delta>0$, there is an admissible $n$-bound $n_0\in\N$.
\end{proposition}

\begin{proof} We demonstrate how to choose $n_0$ so that
$u := P_{n}u_0 = P_{n}A^{-1}b$ solves (rFSM) for every $n\ge n_0$.
For all $m\in\N$ and $n\in\N$, we have
\begin{eqnarray} \nonumber
\|P_{m}AP_{n}u-P_{m}b\|
&=&\|P_{m}AP^2_{n}A^{-1}b-P_{m}b\|\\
\nonumber &\le&\|P_{m}AA^{-1}b-P_{m}b\|
\ +\ \|P_{m}AQ_{n}A^{-1}b\|\\
\nonumber &\le&0\ +\ \|A\|\cdot\|Q_{n}A^{-1}b\|.
\end{eqnarray} But, by our assumption $Q_n\to 0$, there is a
$n_0\in\N$ such that
\begin{equation}
\label{ineq_tau} \|Q_{n}A^{-1}b\|\ \le\ \frac\delta{\|A\|}
\end{equation}
for all $n\ge n_0$, so that $\|P_{m} A P_{n} u - P_{m}b\| < \delta$
holds, and hence $u$ solves (rFSM) for all $n\ge n_0$ and $m\in\N$.
\end{proof}

%=====================================================
%
%=====================================================

\begin{lemma} \label{lem_phi_bound}
Let $n_0\in\N$ be an admissible $n$-bound for $A$, $b$ and a given
precision $\delta>0$. If $n\ge n_0$ and $m\in\N$ are such that
$\|Q_m AP_n\|<1/\|A^{-1}\|$ then the set of all solutions of {\rm
(rFSM)} is a bounded subset of $E$. Precisely, every solution $u\in
\im P_n$ of {\rm (rFSM)} is subject to $\|u\|\le M$ with $M$ given
by (\ref{M}).
\end{lemma}

\begin{proof} Suppose $u\in \im P_n$ solves (rFSM) for given parameters
$\delta,m,n$. Then
\begin{eqnarray}
\nonumber \|Au\| - \|P_m b\| &\le& \|Au - P_m b\|\ =\ \|AP_nu - P_m b\|\\
\nonumber &\le& \|AP_nu - P_m AP_nu\|+\|P_m AP_nu - P_m b\|\\
\nonumber &\le& \|Q_m AP_n\|\cdot\|u\|+\delta
\end{eqnarray}
together with $\|u\|\le\|A^{-1}\|\cdot\|Au\|$ implies that
\begin{eqnarray}
\nonumber \frac{\|u\|}{\|A^{-1}\|}\ \le\ \|Au\|
&\le&\|P_m b\|+\|Q_m AP_n\|\cdot\|u\|+\delta\\
\nonumber &\le&\|b\|+\|Q_m AP_n\|\cdot\|u\|+\delta
\end{eqnarray} and hence
\begin{equation}
\label{M} \|u\|\ \le\ M\ :=\ \frac{\|b\|\ +\ \delta}{1/\|A^{-1}\|\
-\ \|Q_m AP_n\|}.
\end{equation}
\end{proof}

Now we are ready for the key result showing that every solution of
(rFSM) is indeed close to the solution $u_0$ of $Au=b$.

%-----------------------------------------------------
% Convergence Theorem
%-----------------------------------------------------
\begin{theorem} \label{th_conv}
For every $\varepsilon>0$, there are parameters $\delta,m,n$ such
that every solution $u$ of the system {\rm (rFSM)} is an
approximation
\begin{equation}
\label{normapprox} \|u-u_0\|_E<\varepsilon
\end{equation}
of the exact solution $u_0$ of $Au=b$. Precisely, there are three
functions $\delta_0:\R_+\to\R_+$, $n_0:\R_+\to\N$ and
$m_0:\R_+^2\times\N\to\N$ such that if
$\delta<\delta_0(\varepsilon)$, $n\ge n_0(\delta)$ and $m\ge
m_0(\varepsilon,\delta,n)$, then every solution $u\in \im P_n$ of
{\rm (rFSM)} is subject to (\ref{normapprox}).
%\marginpar{Vielleicht etwas zuviel des Guten?}
\end{theorem}

\begin{proof} Let $\varepsilon>0$ be given. We start the proof with three
preliminary steps.
\begin{itemize}
\item[(a)] Choose $\delta<\delta_0:=\frac\varepsilon{3\|A^{-1}\|}$.
\item[(b)] Choose $n_0\in\N$ such that
$\big(\|Q_{n}u_0\|=\big)\|Q_{n}A^{-1}b\|\le\frac\delta{\|A\|}$ for
all $n\ge n_0$, so that $n_0$ is an admissible $n$-bound for
$\delta$ (see inequality (\ref{ineq_tau})). Now let $n\ge n_0$.
\item[(c)] Choose $m_0\in\N$ such that both
$\|Q_{m}b\|<\frac\varepsilon{3\|A^{-1}\|}$ and
\begin{equation}
\label{QAP_bound} \|Q_{m}AP_{n}\|\ <\ \frac
1{\|A^{-1}\|}\left(1-\frac 1{1+\frac
\varepsilon{3(\|b\|+\delta)\cdot\|A^{-1}\|}} \right)
\end{equation}
hold for all $m\ge m_0$, and fix some $m\ge m_0$.
\end{itemize}

Now let $u\in \im P_n$ be a solution of (rFSM) with parameters
$\delta, n$ and $m$ as chosen above. From (\ref{QAP_bound}) we get
$\|Q_{m}AP_{n}\|< 1/\|A^{-1}\|$, and hence, by Lemma
\ref{lem_phi_bound},
\begin{equation}
\label{bound_phi} \|u\|\ \le\ M
\end{equation}
with $M$ as defined in (\ref{M}). Moreover, inequality
(\ref{QAP_bound}) is equivalent to
\[
\|Q_{m}AP_{n}\|\ <\ \frac 1{\|A^{-1}\|}\cdot\frac{\frac
\varepsilon{3(\|b\|+\delta)\cdot\|A^{-1}\|}}{1+\frac
\varepsilon{3(\|b\|+\delta)\cdot\|A^{-1}\|}},
\]
and hence to
\[
\left( {1+\frac
\varepsilon{3(\|b\|+\delta)\cdot\|A^{-1}\|}}\right)\cdot\|Q_{m}AP_{n}\|\
<\ \frac 1{\|A^{-1}\|}\cdot\frac
\varepsilon{3(\|b\|+\delta)\cdot\|A^{-1}\|}.
\]
This, moreover, is equivalent to
\begin{eqnarray}
\nonumber \|Q_{m}AP_{n}\|&<&\frac 1{\|A^{-1}\|}\cdot\frac \varepsilon{3(\|b\|+\delta)\cdot\|A^{-1}\|}-\frac \varepsilon{3(\|b\|+\delta)\cdot\|A^{-1}\|}\cdot\|Q_{m}AP_{n}\|\\
\label{bound_QAP} &=& \frac {\varepsilon\left(
1/\|A^{-1}\|-\|Q_{m}AP_{n}\|\right)}{3(\|b\|+\delta)\cdot\|A^{-1}\|}\
=\ \frac\varepsilon{3M\|A^{-1}\|} \end{eqnarray} with $M$ as defined
in (\ref{M}). Then we have
\begin{eqnarray}
\nonumber \|u-u_0\| &=&\|P_n
u-u_0\|\ =\ \|A^{-1}AP_n u-A^{-1}b\|
\ \le\ \|A^{-1}\|\cdot \|AP_{n}u-b\|\\
\nonumber &\le&\|A^{-1}\|\cdot\big(\|AP_n u-P_{m}AP_{n}u\|\ +\ \|P_{m}AP_{n}u-P_{m}b\|\ +\ \|P_{m}b-b\|\big)\\
\nonumber &<&\|A^{-1}\|\cdot\big(\|Q_{m}AP_{n}\|\cdot\|u\|\ +\ \delta\ +\ \|Q_{m}b\|\big)\\
\nonumber &<&\frac\varepsilon3\ +\ \frac\varepsilon3\ +\
\frac\varepsilon3\ =\ \varepsilon, \end{eqnarray} using inequalities
(\ref{bound_QAP}) and (\ref{bound_phi}) and the bounds on $\delta$
and $\|Q_{m}b\|$ in the last step. \end{proof}

\begin{remark}
One way to effectively solve the system (rFSM) for given parameters
$m,n$ and $\delta$ is to compute a $u\in \im P_n$ that minimizes the
discrepancy in (\ref{eqth:rFSM}), for example using a gradient
method or, if possible, by directly applying the Moore-Penrose
pseudo-inverse $B^+$ of $B:=P_m AP_n$ to the right-hand side $P_m
b$.

If $E$ is a Hilbert space (i.e. if $p=2$ and $X$ is a Hilbert space)
then it is well-known that $u\in \im P_n$ minimizes the residual
$\|Bu - P_m b\|$ if and only if $B^*(Bu - P_m b)=0$. If, in
addition, $P_m$ is self-adjoint for all $m\in\N$, then, after
re-substituting $B$, the latter is equivalent to
\begin{equation}
\label{eq_perturb} P_n A^* P_m A P_n u \ =\ P_n A^* P_m b.
\end{equation}
However, if $m$ is sufficiently large, then, by (\ref{eq:MSMdecay})
and $P_n\to I$, the equation (\ref{eq_perturb}) is just a small
perturbation of
\begin{equation}
\label{eq_fsmAsA} P_n A^* A P_n u \ =\ P_n A^* b,
\end{equation}
which is nothing but the finite section method for the equation
\begin{equation}
\label{eq_AsA} A^* Au\ =\ A^* b.
\end{equation}
Note that the finite section method (\ref{eq_fsmAsA}) is applicable
since $A^* A$ is positive definite (see, e.g. Theorem 1.10 b in
\cite{HaRoSi:2001}). Clearly, if $A$ is invertible, as we require,
then also its adjoint $A^*$ is invertible, and (\ref{eq_AsA}) is
equivalent to our original equation $Au=b$.

Summarizing, if $E$ is a Hilbert space and all $P_m$ are
self-adjoint, then minimizing $\|P_m A P_n u - P_m b\|$ is
equivalent to solving a slight perturbation (\ref{eq_perturb}) of
the finite section method (\ref{eq_fsmAsA}) for (\ref{eq_AsA}).
\end{remark}

%%%%%%%%%%%%%%%%%%%%%%%%%%%%%%%%%%%%%%%%%%%%%%%%%%%%%%%%%%%%%%%%%%%
\section{Summary}
Compared to the FSM and its subsequence version fully characterized
in Theorem \ref{th:main}, the rFSM imposes no further conditions on
the operator $A$ (such as richness, conditions on its limit
operators, etc.) other than its invertibility and the rather mild
decay property (\ref{eq:MSMdecay}). Of course, on the down side, we
are restricted to $p<\infty$, and, more seriously, for general
operators $A$ we do not really know yet how to choose $m$ in
dependence on $n$. However, the choice $m=n+w$ is clear for
operators with band-width $w$, and something similar is possible for
a band-dominated operator $A$ (where $w$ must be fitted to the
function $f_A$ from \cite[p. 32ff]{Li:Book}). In contrast, in
\cite{HeiLiPott:MSM} we have chosen $m=\frac 65 n$.

We close with the following example.

\begin{example}
Let $E=\ell^2(\Z,\C)$, $b=(b(i))_{i\in\Z}\in E$ with
$b(i)=2^{-|i|}$, and let $A:E\to E$ be given in block matrix
notation by
\[
A\ =\ \left(\begin{array}{ccccc}
\ddots&\ddots\\
&B&C\\
&&B&C\\
&&&B&\ddots\\
&&&&\ddots
\end{array}\right),
\]
where
\[
B\ =\ \left(\begin{array}{ccc} 1&1&0\\1&0&0\\0&0&0\end{array}\right)
\qquad\textrm{and}\qquad C\ =\ \left(\begin{array}{ccc}
0&0&0\\0&0&0\\1&1&1\end{array}\right),
\]
and where one of the $B$ blocks is located at position
$\{(i,j):i,j\in\{-1,0,1\}\}$ in $A$. Then $A$ is a band operator
with band-width $w=3$, it is invertible but has a nonzero plus-index
$\ind_+(A)=\ind(P_\N AP_\N|_\im P_\N)=1$. From Proposition
\ref{prop:N=1:I} we know that there is no choice of $\Omega$ and
$\I\subset\N$ that makes the FSM (\ref{eq:FSsubs}) for $Au=b$
stable. But Remark \ref{rem:Gohberg} shows that passing to the
equivalent system $V_1Au=V_1b$ might solve this problem since
$A':=V_1 A$ has plus-index zero. This approach leads us to looking
at finite sections of
\[
A'\ =\ V_1A\ =\ \left(\begin{array}{ccccc}
\ddots\\
&D\\
&&D\\
&&&D\\
&&&&\ddots
\end{array}\right)
\qquad\textrm{with}\qquad D\ =\ \left(\begin{array}{ccc}
1&1&1\\1&1&0\\1&0&0\end{array}\right)
\]
and one of the $D$ blocks located at position
$\{(i,j):i,j\in\{-1,0,1\}\}$ in $A'$. For $\Omega=[-1,1]$, we get
that $\sstab_{\Omega,\N}(A')=\{A',F,G,H,J,K,L\}$, where
\[
\begin{array}{rl}
F =
\diag\left(\cdots,1,\underline{1},\left(\begin{array}{cc}1&0\\0&0\end{array}\right),D,D,\cdots\right),
& J =
\diag\left(\cdots,D,D,\left(\begin{array}{cc}1&1\\1&1\end{array}\right),\overline{1},1,\cdots\right),
\\
G = \diag\left(\cdots,1,\underline{1},D,D,D,\cdots\right),&
K = \diag\left(\cdots,D,D,D,\overline{1},1,\cdots\right), \\
H = \diag\left(\cdots,1,\underline{1},0,D,D,\cdots\right), & L =
\diag\left(\cdots,D,D,1,\overline{1},1,\cdots\right)
\end{array}
\]
with the underlined $1$'s at position $(-1,-1)$ and the overlined
$1$'s at position $(1,1)$ in the respective matrix. Out of the seven
elements of $\sstab_{\Omega,\N}(A')$ only $A',G,K$ and $L$ are
invertible (noting that $D$ is invertible). Looking at certain index
subsets $\I$ of $\N$, we see that
\[
\sstab_{\Omega,3\N}(A')\, =\, \{A',F,J\},\qquad
\sstab_{\Omega,3\N+1}(A')\, =\, \{A',G,K\},\qquad
\sstab_{\Omega,3\N+2}(A')\, =\, \{A',H,L\},
\]
so that the set $\I=3\N+1$ yields a stable subsequence
\eqref{eq:FSsubs}, whereas $3\N$ and $3\N+2$ don't. After these
considerations (shift the system by 1 row, single out a stable
subsequence) it is now straightforward to approximately solve the
equation via the FSM.

In contrast, the rFSM immediately applies to $Au=b$ if we choose
$m=n+3$ (since $A$ has band-width $w=3$) and solve the systems
$P_{n+3}AP_n u\approx P_{n+3}b$ approximately for $n=1,2,...$ using
the Moore-Penrose pseudo-inverse of $P_{n+3}AP_n$. The latter means
that we compute
\[
u_n\ :=\ (P_{n+3}AP_n)^+\, P_{n+3}b,\qquad n=1,2,...
\]
and get that $u_n$ (if extended by zeros to an infinite vector)
converges to the exact solution $u_0=A^{-1}b$ with
\begin{eqnarray*}
\|u_n-u_0\| &\le& \|A^{-1}\|\
(\|P_{n+3}AP_nu_n-P_{n+3}b\|+\|Q_{n+3}b\|)\\
&\le& 2\,\left(\frac {3\cdot 2}{2^{n}}\ +\ \frac 1{2^{n+3}}\right)\
=\ \frac{49}{2^{n+2}},
\end{eqnarray*}
which follows from the computations in the proofs of Proposition
\ref{prop:exist} and Theorem \ref{th_conv} together with $\|A\|\le
3$ and $\|A^{-1}\|\le 2$.
\end{example}

\bigskip

{\bf Acknowledgements. } I would like to thank the organizers,
especially Robert Edward and Golden Thambithurai, for their
invitation to and great hospitality at the International Conference
on Functional Analysis at Nagercoil, India. The research presented
in this paper was financially supported by Marie-Curie Grants
MEIF-CT-2005-009758 and PERG02-GA-2007-224761 of the EU and was
largely carried out on the very pleasant trip to Nagercoil.
%%%%%%%%%%%%%%%%%%%%%%%%%%%%%%%%%%%%%%%%%%%%%%%%%%%%%%%%%%%%%%%%%%%%%%%%%%%%

\bigskip

\noindent {\bf Author's address:}\\[2mm]
Marko Lindner\hfill {\tt marko.lindner@mathematik.tu-chemnitz.de}\\
Fakult\"at Mathematik\\
TU Chemnitz\\
D-09107 Chemnitz\\
GERMANY

\end{document}